\documentclass{article}

\usepackage{amsfonts}
\usepackage{amssymb}
\usepackage{hyperref}
\usepackage{amsthm}
\usepackage{latexsym}
\usepackage{xspace}
\usepackage{graphicx}
\usepackage{amsmath}
\usepackage[latin1]{inputenc}
\usepackage{url}
\usepackage{epsfig}

\newtheorem{dfn}{Definition}[section]
\newtheorem{ex}{Example}[section]
\newtheorem{thm}{Theorem}[section]

\newtheorem{cor}{Corollary}[section]

\newcommand{\sd}{\stackrel{d}{=}}

\newcommand{\bx}{\boldsymbol{x}}

\newcommand{\by}{\boldsymbol{y}}
\newcommand{\bz}{\boldsymbol{z}}

\newcommand{\bu}{\boldsymbol{u}}

\newcommand{\bzero}{\boldsymbol{0}}

\begin{document}
\title{Array Variate Elliptical Random Variables with Multiway Kronecker Delta Covariance Matrix Structure \\ Deniz Akdemir \\ \small Department of Statistics \\ 
  University of Central Florida\\ Orlando, FL 32816 \normalsize }

\maketitle

\begin{abstract} Standard statistical methods applied to matrix random variables often fail to describe the underlying structure in multiway data sets. In this paper we will discuss the concept of an array variate random variable and introduce a class of elliptical array densities which have elliptical contours.\end{abstract}

%\amssubj{AMS 2000 Subject Classification: Primary 62H10, Secondary 62H05.}

%\keywords{Keywords \& Phrases: Random Matrices, Random Arrays, Array Variate Normal Distribution, Array Variate Elliptical Distribution, Multilevel Data Analysis, Repeated Measures, Classification, Dimension Reduction}
\allowdisplaybreaks

\section{Introduction}

The array variate random variable up to $2$ dimensions has been studied intensively in \cite{gupta2000matrix} and by many others. For arrays observations of $3,$ $4$ or in general $i$ dimensions a suitable normal probablity model has been recently proposed in \cite{DenizGupta}. In this paper we will generalize the notion of elliptical random variables to the array case.

In Section 2, we first study the algebra of arrays. In Section 3, we introduce the concept of an array variable random variable. In section 4, the density of the elliptical array random variable is provided. Then, in Section 5, we define the array random variable with elliptical density.

\section{Array Algebra}

In this paper we will only study arrays with real elements. We will write $\widetilde{X}$ to say that $\widetilde{X}$ is an array. When it is necessary we can write the dimensions of the array as subindices, e.g., if $\widetilde{X}$ is a $m_1 \times m_2\times m_3 \times m_4$ dimensional array in $R^{m_1\times m_2 \times \ldots \times m_i}$, then we can write $\widetilde{X}_{m_1 \times m_2\times m_3 \times m_4}.$  Arrays with the usual elementwise summation and scalar multiplication operations can be shown to be a vector space.  

To refer to an element of an array $\widetilde{X}_{m_1 \times m_2\times m_3 \times m_4},$ we write the position of the element as a subindex to the array name in parenthesis, $(\widetilde{X})_{r_1r_2r_3r_4}.$  If we want to refer to a specific column vector obtained by keeping all but an indicated dimension constant, we indicate the constant dimensions as before but we will put '$:$' for the non constant dimension, e.g., for  $\widetilde{X}_{m_1 \times m_2\times m_3 \times m_4},$    $(\widetilde{X})_{r_1r_2:r_4}$ refers to the the column vector $(({X})_{r_1r_21r_4},({X})_{r_1r_22r_4}, \ldots, ({X})_{r_1r_2m_3r_4})'.$
 
We will now review some basic principles and techniques of array algebra. These results and their proofs can be found in Rauhala \cite{rauhala1974array},  \cite{rauhala1980introduction}  and Blaha \cite{blaha1977few}. 

\begin{dfn}\emph{Inverse Kronecker product} of two matrices $A$ and $B$ of dimensions $p\times q$ and $r \times s$ correspondingly is written as $A\otimes^i B$ and is defined as $A\otimes^i B=[A(B)_{jk}]_{pr\times qs}=B\otimes A,$ where $'\otimes'$ represents the ordinary Kronecker product.\end{dfn}

The following properties of the inverse Kronecker product are useful:\begin{itemize}\item $\bzero \otimes^i A= A \otimes^i \bzero=\bzero.$ \item $(A_1+A_2)\otimes^i B=A_1\otimes^i B+ A_2 \otimes^i B.$ \item $A \otimes^i (B_1+ B_2)=A \otimes^i B_1+ A \otimes^i B_2.$ \item  $\alpha A \otimes^i \beta B= \alpha \beta A \otimes^i B.$ \item $(A_1 \otimes^i B_1)(A_2 \otimes^i B_2)= A_1A_2 \otimes^i B_1B_2.$ \item  $(A \otimes^i B)^{-1}=(A^{-1} \otimes^i B^{-1}).$ \item  $(A \otimes^i B)^{+}=(A^{+} \otimes^i B^{+}),$ where $A^{+}$ is the Moore-Penrose inverse of $A.$ \item $(A \otimes^i B)^{-}=(A^{-} \otimes^i B^{-}),$ where $A^{-}$ is the $l$-inverse of $A$ defined as $A^{-}=(A'A)^{-1}A'.$\item If $\{\lambda_i\}$ and $\{\mu_j\}$ are the eigenvalues with the corresponding eigenvectors $\{\bx_i\}$ and $\{\by_j\}$ for matrices $A$ and $B$ respectively, then $A\otimes^i B$ has eigenvalues $\{\lambda_i\mu_j\}$ with corresponding eigenvectors $\{\bx_i\otimes^i\by_j\}.$\item Given two matrices $A_{n\times n}$ and $B_{m\times m}$ $|A\otimes^i B|=|A|^m|B|^n,$ $tr(A\otimes^i B)=tr(A)tr(B).$\item $A\otimes^i B=B\otimes A=U_1 A\otimes B U_2,$ for some permutation matrices $U_1$ and $U_2.$ 
\end{itemize}

It is well known that a matrix equation $$AXB'=C$$ can be rewritten in its monolinear form as \begin{equation}\label{eqmnf}A\otimes^i B vec(X)=vec(C).\end{equation} Furthermore, the matrix equality $$A\otimes^i B XC'=E$$ obtained by stacking equations of the form (\ref{eqmnf}) can be written in its monolinear form as $$(A\otimes^i B \otimes^i C) vec(X)=vec(E).$$ This process of stacking equations could be continued and R-matrix multiplication operation introduced by Rauhala \cite{rauhala1974array} provides a compact way of representing these equations in array form:

\begin{dfn}\emph{R-Matrix Multiplication} is defined elementwise: 

 $$((A_1)^1 (A_2)^2 \ldots (A_i)^i\widetilde{X}_{m_1 \times m_2 \times \ldots \times m_i})_{q_1q_2\ldots q_i}$$ $$=\sum_{r_1=1}^{m_1}(A_1)_{q_1r_1}\sum_{r_2=1}^{m_2}(A_2)_{q_2r_2}\sum_{r_3=1}^{m_3}(A_3)_{q_3r_3}\ldots \sum_{r_i=1}^{m_i}(A_i)_{q_ir_i}(\widetilde{X})_{r_1r_2\ldots r_i}.$$ \end{dfn}

R-Matrix multiplication generalizes the matrix multiplication (array multiplication in two dimensions)to the case of $k$-dimensional arrays. The following useful properties of the R-Matrix multiplication are reviewed by Blaha \cite{blaha1977few}:
\begin{enumerate}
\item $(A)^1B=AB.$ 
\item $(A_1)^1(A_2)^2C=A_1CA'_2.$ 
\item $\widetilde{Y}=(I)^1(I)^2\ldots (I)^i \widetilde{Y}.$ 
\item \small $((A_1)^1 (A_2)^2\ldots (A_i)^i)((B_1)^1(B_2)^2\ldots (B_i)^i)\widetilde{Y}= (A_1B_1)^1(A_2B_2)^2\ldots(A_iB_i)^i\widetilde{Y}.$ \normalsize
\end{enumerate}

The operator $rvec$ describes the relationship between $\widetilde{X}_{m_1 \times m_2 \times \ldots m_i}$ and its monolinear form $\bx_{m_1m_2\ldots m_i\times 1}.$ 
\begin{dfn}\label{def:rvec} $rvec( \widetilde{X}_{m_1 \times m_2 \times \ldots m_i})=\bx_{m_1m_2\ldots m_i\times 1}$ where $\bx$ is the column vector obtained by stacking the elements of the array $\widetilde{X}$ in the order of its dimensions; i.e., $(\widetilde{X})_{j_1 j_2 \ldots j_i}=(\bx)_j$ where $j=(j_i-1)n_{i-1}n_{i-2}\ldots n_1+(j_i-2)n_{i-2}n_{i-3}\ldots n_1+\ldots+(j_2-1)n_1+j_1.$\end{dfn}

\begin{thm}\label{rmultvec}Let $\widetilde{L}_{m_1 \times m_2 \times\ldots m_i}=(A_1)^1(A_2)^2\ldots(A_i)^i\widetilde{X}$ where $(A_j)^j$ is an $m_j\times n_j$ matrix for $j=1,2,\ldots,i$ and $\widetilde{X}$ is an $n_1\times n_2\times\ldots\times n_i$ array. Write $\mathbf{l}=rvec(\widetilde{L})$ and $\bx=rvec(\widetilde{X}).$ Then, $\mathbf{l}=A_1\otimes^iA_2\otimes^i\ldots\otimes^i A_i\bx.$ \end{thm}

Therefore, there is an equivalent expression of the array equation in monolinear form.

\begin{dfn}{} The square norm of $\widetilde{X}_{m_1 \times m_2 \times\ldots m_i}$ is defined as $$\|\widetilde{X}\|^2=\sum_{j_1=1}^{m_1}\sum_{j_2=1}^{m_2}\ldots\sum_{j_i=1}^{m_i}((\widetilde{X})_{j_1j_2\ldots j_i})^2.$$ \end{dfn}

\begin{dfn}{} The distance of $\widetilde{X_1}_{m_1 \times m_2 \times\ldots m_i}$ from $\widetilde{X_2}_{m_1 \times m_2 \times\ldots m_i}$ is defined as $\sqrt{\|\widetilde{X_1}-\widetilde{X_2}\|^2}.$ \end{dfn}

\begin{ex} Let $\widetilde{Y}=(A_1)^1 (A_2)^2\ldots (A_i)^i\widetilde{X}+\widetilde{E}.$ Then $\|\widetilde{E}\|^2$ is minimized for $\widehat{\widetilde{X}}=(A_1^{-})^1(A_2^{-})^2\ldots(A_i^{-})^i\widetilde{Y}.$ \end{ex}
 
\section{Array Variate Random Variables}

Arrays can be constant arrays, i.e. if $(\widetilde{X})_{r_1r_2\ldots r_i}\in\mathbf{R}$ are constants for all $r_j,$ $j=1,2,\ldots,m_j$ and $j=1,2,\ldots, i$ then the array $\widetilde{X}$ is a constant array. 

Array variate random variables are arrays with all elements $(\widetilde{X})_{r_1r_2\ldots r_i}\in\mathbf{R}$ random variables. If the sample space for the random outcome $s$ is $\mathbb{S}$, $(\widetilde{X})_{r_1r_2\ldots r_i}=(\widetilde{X}(s))_{r_1r_2\ldots r_i}$ where each of $(\widetilde{X}(s))_{r_1r_2\ldots r_i}$ is a real valued function from $\mathbb{S}$ to $\mathbb{R}.$  

If $\widetilde{X}$ is an array variate random variable, its density (if it exists) is a scalar function $f_{\widetilde{X}}(\widetilde{X})$ such that:

\begin{itemize}
\item $f_{\widetilde{X}}(\widetilde{X})\geq 0;$ 
\item $\int_{\widetilde{X}} f_{\widetilde{X}}(\widetilde{X})d\widetilde{X}=1;$ 
\item $P(\widetilde{X}\in A)=\int_{A} f_{\widetilde{X}}(\widetilde{X})d\widetilde{X},$ where A is a subset of the space of realizations for $\widetilde{X}.$ 
\end{itemize}
 
A scalar function $f_{\widetilde{X},\widetilde{Y}}(\widetilde{X},\widetilde{Y})$ defines a joint (bi-array variate) probability density function if 
 \begin{itemize}
\item $f_{\widetilde{X},\widetilde{Y}}(\widetilde{X},\widetilde{Y})\geq 0;$ 
\item $\int_{\widetilde{Y}}\int_{\widetilde{X}} f_{\widetilde{X},\widetilde{Y}}(\widetilde{X},\widetilde{Y})d\widetilde{X}d\widetilde{Y}=1;$ 
\item $P((\widetilde{X},\widetilde{Y})\in A)=\int \int_{A} f_{\widetilde{X},\widetilde{Y}}(\widetilde{X},\widetilde{Y})d\widetilde{X}d\widetilde{Y},$ where A is a subset of the space of realizations for $(\widetilde{X},\widetilde{Y}).$ 
\end{itemize}
The marginal probability density function of $\widetilde{X}$ is defined by $$f_{\widetilde{X}}(\widetilde{X})=\int_{\widetilde{Y}} f_{\widetilde{X},\widetilde{Y}}(\widetilde{X},\widetilde{Y})d\widetilde{Y},$$ and the conditional probability density function of $\widetilde{X}$ given $\widetilde{Y}$ is defined by $$f_{\widetilde{X}|\widetilde{Y}}(\widetilde{X}|\widetilde{Y})=\frac{f_{\widetilde{X},\widetilde{Y}}(\widetilde{X},\widetilde{Y})}{f_{\widetilde{Y}}(\widetilde{Y})},$$ where $f_{\widetilde{Y}}(\widetilde{Y})>0.$ 

Two random arrays $\widetilde{X}$ and  $\widetilde{Y}$ are independent if and only if $$f_{\widetilde{X},\widetilde{Y}}(\widetilde{X},\widetilde{Y})= f_{\widetilde{X}}(\widetilde{X})f_{\widetilde{Y}}(\widetilde{Y}).$$

\begin{thm}{}  Let $(A_1)^1,$ $(A_2)^2,$ $\ldots,$ $(A_i)^i$ be $m_1,$ $m_2,$ $\ldots,$ $m_i$ dimensional positive definite matrices. The Jacobian $J(\widetilde{X}\rightarrow \widetilde{Z})$ of the transformation $\widetilde{X}=(A_1)^1 (A_2)^2 \ldots (A_i)^i\widetilde{Z}+\widetilde{M}$ is  $$(|A_1|^{\prod_{j\neq 1}{m_j}} |A_2|^{\prod_{j\neq 2}{m_j}} \ldots |A_i|^{\prod_{j\neq i}{m_j}})^{-1}.$$
\begin{proof} The result is proven using the equivalence of monolinear form obtained through the $rvec(\widetilde{X})$ and array $\widetilde{X}.$ Let $\widetilde{L}_{m_1 \times m_2 \times\ldots m_i}=(A_1)^1(A_2)^2\ldots(A_i)^i\widetilde{Z}$ where $(A_j)^j$ is an $m_j\times n_j$ matrix for $j=1,2,\ldots,i$ and $\widetilde{X}$ is an $n_1\times n_2\times\ldots\times n_i$ array. Write $\mathbf{l}=rvec(\widetilde{L})$ and $\bz=rvec(\widetilde{Z}).$ Then, $\mathbf{l}=A_1\otimes^iA_2\otimes^i\ldots\otimes^i A_i\bz.$ The result follows from noting that $J(\mathbf{l}\rightarrow \bz)=|A_1\otimes^iA_2\otimes^i\ldots\otimes^i A_i|^{-1},$ and using induction with the rule $|A\otimes^i B|=|A|^m|B|^n$ for ${n\times n}$ matrix $A$ and ${m\times m}$ matrix $B$ to show that \small $|A_1\otimes^iA_2\otimes^i\ldots\otimes^i A_i|^{-1}=(|A_1|^{\prod_{j\neq 1}{m_j}} |A_2|^{\prod_{j\neq 2}{m_j}} \ldots |A_i|^{\prod_{j\neq i}{m_j}})^{-1}.$ \normalsize \end{proof}\end{thm}

\begin{cor}{\label{jacobian}} Let $\widetilde{Z}\sim f_{\widetilde{Z}} (\widetilde{Z}).$ Define $\widetilde{X}=(A_1)^1 (A_2)^2 \ldots (A_i)^i\widetilde{Z}+\widetilde{M}$ where $(A_1)^1, (A_2)^2, \ldots, (A_i)^i$ be $m_1, m_2, \ldots, m_i$ dimensional positive definite matrices. The pdf of $\widetilde{X}$ is given by \small $$f_{\widetilde{X}} (\widetilde{X}; (A_1)^1, (A_2)^2, \ldots, (A_i)^i, \widetilde{M})=\frac{f(A_1^{-1})^1 (A_2^{-1})^2 \ldots (A_i^{-1})^i(\widetilde{X}-\widetilde{M}))}{|A_1|^{\prod_{j\neq 1}{m_j}} |A_2|^{\prod_{j\neq 2}{m_j}} \ldots |A_i|^{\prod_{j\neq i}{m_j}}} .$$ \normalsize \end{cor}

The main advantage in choosing a Kronecker structure is the decrease in the number of parameters.

\section{Array Variate Normal Distribution}

By using the results in the previous section on array algebra, mainly the relationship of the arrays to their monolinear forms described by Definition \ref{def:rvec} , we can write the density of the standard normal array variable.
\begin{dfn}{} If $$\widetilde{Z}\sim N_{m_1 \times m_2\times \ldots \times m_i}(\widetilde{M}=\bzero, \Lambda=I_{m_1m_2\ldots m_i}),$$ then $\widetilde{Z}$ has array variate standard normal distribution. The pdf of $\widetilde{Z}$ is given by \begin{equation}\label{eq:densitystarn}f_{\widetilde{Z}} (\widetilde{Z})=\frac{\exp{(-\frac{1}{2}\|\widetilde{Z}\|^2)}}{(2\pi)^{m_1m_2\ldots m_i/2}}.\end{equation}\end{dfn}

For the scalar case, the density for the standard normal variable $z \in \mathbf{R}^1$ is given as \[\phi_1(z)=\frac{1}{(2\pi)^\frac{1}{2}}exp(-\frac{1}{2}z^2).\] For the $m_1$ dimensional  standard normal vector $\bz\in\mathbf{R}^{m_1},$ the density is given by \[\phi_{m_1}(\bz)=\frac{1}{(2\pi)^\frac{m_1}{2}}exp(-\frac{1}{2}\bz'\bz).\] Finally the $m_1\times m_2$ standard matrix variate variable $Z\in\mathbf{R}^{m_1\times m_2}$ has the density \[\phi_{m_1\times m_2}(Z)=\frac{1}{(2\pi)^\frac{m_1m_2}{2}}exp(-\frac{1}{2}trace(Z'Z)).\] With the above definition, we have generalized the notion of normal random variable to the array variate case.

\begin{dfn}{\label{stdnorm}} We write $$\widetilde{X}\sim N_{m_1 \times m_2\times \ldots \times m_i}(\widetilde{M}, \Lambda_{m_1m_2\ldots m_i})$$ if   $rvec(\widetilde{X})\sim N_{m_1m_2\ldots m_i}(rvec(\widetilde{M}),\Lambda_{m_1m_2\ldots m_i}).$ Here, $\widetilde{M}$ is the expected value of $\widetilde{X}$, and $\Lambda_{m_1m_2\ldots m_i}$ is the covariance matrix of the ${m_1m_2\ldots m_i}$-variate random variable $rvec(\widetilde{X}).$\end{dfn} 

The family of normal densities with Kronecker Delta Covariance Structure are obtained by considering the densities obtained by the location-scale transformations of the standard normal variables. This kind of model is defined in the next.

\begin{thm}{}\label{modkroncov} Let $\widetilde{Z}\sim N_{m_1 \times m_2\times \ldots \times m_i}(\widetilde{M}=\bzero, \Lambda=I_{m_1m_2\ldots m_i}).$ Define $\widetilde{X}=(A_1)^1 (A_2)^2 \ldots (A_i)^i\widetilde{Z}+\widetilde{M}$ where $A_1, A_2,\ldots,A_i$ are non singular matrices of orders $m_1, m_2,\ldots, m_i$.   Then the pdf of $\widetilde{X}$ is given by \small \begin{equation}\label{eq:densityarn}\phi(\widetilde{X}; \widetilde{M},A_1,A_2,\ldots A_i)=\frac{\exp{(-\frac{1}{2}\|{(A_1^{-1})^1 (A_2^{-1})^2 \ldots (A_i^{-1})^i(\widetilde{X}-\widetilde{M})}\|^2)}}{(2\pi)^{m_1m_2\ldots m_i/2}|A_1|^{\prod_{j\neq 1}{m_j}} |A_2|^{\prod_{j\neq 2}{m_j}} \ldots |A_i|^{\prod_{j\neq i}{m_j}}}.\end{equation} \normalsize \begin{proof} The result is easily obtained using the density in Definition \ref{stdnorm} with Corollary \ref{jacobian}.\end{proof}\end{thm}

\section{Array Variate Elliptical Distribution}

A spherically symmetric random vector $\bx$ has a density of the form $f(\bx'\bx)$ for some kernel pdf $f(t),$ defined for $t\in \mathbf{R}^+.$ 

A stochastic representation for the random vector $\bx$ is given by $\bx\sd r\bu,$ where $r$ is a nonnegative random variable with cdf $K(r)$ that is independent of  $\bu$ which is uniformly distributed over the unit sphere in $\mathbf{R}^k,$ denoted by $\mathbf{S}^{k},$ where  $r\sd \sqrt(\bx'\bx),$ $\bu\sd\frac{\bx}{\sqrt(\bx'\bx)}.$ If both $K'(r)=k(r)\in\mathbf{K}$ the pdf of $r,$ and  $f(\bx'\bx)$ the pdf of $\bx$ exist, then they are related as follows: \begin{equation}k(r)=\frac{2\pi^{\frac{k}{2}}}{\Gamma(\frac{k}{2})}r^{k-1}f(r^2).\end{equation}

Some examples follow:
\begin{enumerate}
	\item The standard multivariate normal $N_k(\bzero_k, I_k)$ distribution with pdf $$\phi_k(\bx)=\frac{1}{(2\pi)^{k/2}}e^{-\frac{1}{2}\bx'\bx}.$$ 
 	\item The multivariate spherical $t$ distribution with $v$ degrees of freedom with density function $$f(\bx)=\frac{\Gamma(\frac{1}{2}(v+k))}{\Gamma(\frac{1}{2}v)(v\pi)^{k/2}}\frac{1}{(1+\frac{1}{v}\bx'\bx)^{(v+k)/2}}.$$
 \item When $v=1$ the multivariate $t$ distribution is called spherical Cauchy distribution.
\end{enumerate}

The array variate random variable with spherical contours is defined using the relationship between the array and its monolinear form. This amounts to saying that the density of the array variable $\widetilde{X}$ is spherical if and only if it can be written in the form $f(\|\widetilde{X}\|^2)$ for some kernel pdf $f(t),$ defined for $t\in \mathbf{R}^+.$ The elliptically contoured random variables with unstructured covariance matrices can be obtained by applying a linear transformation to the monolinear form of an array variate spherical random variable. 

In general we can define an array variable random variable with unstructured covariance matrix, i.e. we can say $\widetilde{X}$ has array variate distribution if $rvec(\widetilde{X})$ has  elliptical distribution. Array variate densities with elliptical contours that have Kronecker delta covariance structure are also easily constructed using Corrolary \ref{jacobian}. 

\begin{dfn}\label{elliptical} Let  $A_1, A_2,\ldots,A_i$ be non singular matrices with orders $m_1,$  $m_2,$ $\ldots,$ $m_i$ and $\widetilde{M}$ be a $m_1 \times$ $m_2$ $\times \ldots$ $\times m_i$ constant array. Also, let $m=m_1m_2\ldots m_i.$ Then the pdf of an $m_1 \times$ $m_2$ $\times \ldots$ $\times m_i$ elliptically contoured array variate random variable $\widetilde{X}$ with kernel  pdf $f(t),$ $t>0$ is  \small $$f_{\widetilde{X}} (\widetilde{X}; (A_1)^1, (A_2)^2, \ldots, (A_i)^i, \widetilde{M})=\frac{f(\|(A_1^{-1})^1 (A_2^{-1})^2 \ldots (A_i^{-1})^i(\widetilde{X}-\widetilde{M})\|^2)}{|A_1|^{\prod_{j\neq 1}{m_j}} |A_2|^{\prod_{j\neq 2}{m_j}} \ldots |A_i|^{\prod_{j\neq i}{m_j}}} .$$ \normalsize   
\end{dfn}

For example, the following definition provides a generalization of the multivariate $t$ distribution to the array variate case.

\begin{dfn} Let  $A_1, A_2,\ldots,A_i$ be non singular matrices with orders $m_1,$  $m_2,$ $\ldots,$ $m_i$ and $\widetilde{M}$ be a $m_1 \times$ $m_2$ $\times \ldots$ $\times m_i$ constant array. Also, let $m=m_1m_2\ldots m_i.$ Then the pdf of an $m_1 \times$ $m_2$ $\times \ldots$ $\times m_i$ array variate t random variable $\widetilde{T}$ with degrees of freedom $v$ is given by \small \begin{equation}\label{eq:densityart}f(\widetilde{T}; \widetilde{M},A_1,A_2,\ldots A_i)=c\frac{(1+\|{(A_1^{-1})^1 (A_2^{-1})^2 \ldots (A_i^{-1})^i(\widetilde{T}-\widetilde{M})\|^2)^{-(v+m)/2}}}{|A_1|^{\prod_{j\neq 1}{m_j}} |A_2|^{\prod_{j\neq 2}{m_j}} \ldots |A_i|^{\prod_{j\neq i}{m_j}}}\end{equation}  \normalsize where $c=\frac{(v\pi)^{m/2}\Gamma((v+m)/2)}{\Gamma(v/2)}.$  
\end{dfn}

Distributional properties of a array normal variable with density in the form of given by Theorem \ref{modkroncov} or a array variate elliptical random variable  with density of the form in (\ref{elliptical}) can obtained by using the equivalent monolinear representation of the random variable. The moments, the marginal and conditional distributions, independence of variates can be studied considering the equivalent monolinear form of the array variable and the well known properties of the multivariate normal random variable and the multivariate elliptical random variable. 
\bibliographystyle{plain}
\bibliography{arrayref}
\end{document}